\documentclass[a4paper,10pt]{article} 
\usepackage[english]{babel} 
\usepackage{amsmath}
\usepackage{amssymb,amsthm}
\usepackage{enumerate}

\newtheorem{thm}{Theorem}
\newtheorem{lemma}[thm]{Lemma}
\newtheorem{prop}[thm]{Proposition}
\newtheorem{cor}[thm]{Corollary}

\newtheorem*{conj}{Conjecture}

\theoremstyle{definition}
\newtheorem{defi}[thm]{Definition}

\theoremstyle{remark}
\newtheorem*{remark}{Remark}

\newtheorem{example}[subsubsection]{Example}

\newcommand{\Z}{\mathbb{Z}}
\newcommand{\Q}{\mathbb{Q}}
\newcommand{\C}{\mathbb{C}}

\def\mapdown#1{\Big\downarrow \rlap{$\vcenter{\hbox{$\scriptstyle#1$}}$}}

\newcommand{\da}{\dashrightarrow}
\newcommand{\la}{\longrightarrow}

\def\Mg{\overline{M_g}}

\def\CF{{\mathcal C}_g(F)}
\def\CL{{\mathcal C}_g(L)}
\def\CFS{{\mathcal C}_g(F,S)}
\def\ML{{\mathcal M}_g(L)}
\def\MLS{{\mathcal M}_g(L,S)}
\def\CLS{{\mathcal C}_g(L,S)}

\def\O{\mathcal O}

\def\ni{\noindent}
\def\cl{\centerline}

\def\ub{B_g(K)}
\def\gb{N_g}

\def\C{{\Bbb C}}
\def\Q{{\Bbb Q}}
\def\P{{\Bbb P}}

\def\Z{{\Bbb Z}}

\def\O{{\cal O}}

\def\Im{\mathop{\rm Im}}

\def\Spec{\mathop{\rm Spec }}

\def\Aut{\mathop{\rm Aut }}

\def\Xn{X^n_B}
\def\bir{\approx }

\begin{document}
\cl{{\bf MODULI THEORY AND ARITHMETIC OF ALGEBRAIC VARIETIES}\footnote{To appear in the
Proceedings of the Fano Conference}}

\

\cl{Lucia Caporaso}

\

\tableofcontents 

\

\

\section{Introduction and Preliminaries}

\subsection{Summary}
This paper surveys a few applications of algebro-geometric moduli theory
to issues concerning the distribution of rational points in algebraic varieties.

A few well known  arithmetic problems with their expected answers 
(the so-called ``diophantine conjectures") are introduced in  section~\ref{dio},
explaining their connection with  a  circle of ideas, 
whose goal is   to find a unifying
theme
in analytic, differential, algebraic  and arithmetic geometry.

In section~\ref{AG} we illustrate how, applying the modern moduli theory of  algebraic
varieties,
 the diophantine conjectures imply  uniform boundedness of rational points on curves,
which is also  open. The crucial tool is an algebro-geometric theorem
(the ``Correlation" or ``fibered power" theorem) about families of varieties of general type.

A series of enumerative results  about rational points on curves is described in  section~\ref{count}.

Finally, we examine related  issues over  function fields.
This  is  an appropriate thing to do, in fact, as we explain in section
\ref{ffdio}, part of the motivation for the diophantine conjectures  
came from parallel results over function fields.

We show how  other aspects of moduli theory
 are used to obtain  
that various uniform boundedness phenomena do hold over function fields.

\subsection{Conventions}
\label{prel}
We work in characteristic zero.
 The following notation will be fixed throughout. $K$ is a number field and 
 $V$ 
an algebraic variety defined over
$K$; $V$ is projective, reduced,
irreducible, and of positive dimension. The set of
$K$-rational points of
$V$ is denoted by
$V(K)$.

When considering $V$ over an extension of $K$, for example over $\C$, we will continue to denote it by
$V$.

All varieties are assumed to be connected. The word ``curve" means  projective variety of 
dimension 1.
The word ``family"  stands for a projective morphism of integral  varieties,
often denoted by
$f: X\la B$; the fiber over  a point $b\in B$ is denoted by $X_b:= f^{-1}(b)$.

A complex function field,    denoted by $L=\C(B)$, is the field
of rational functions of an integral variety $B$, which can be assumed to
be projective, to fix ideas. 

Whenever dealing with  number fields   and function fields  simultaneoulsly,
we will use the special notation $F$ to denote either type of field.

The symbols $g,\ d,\ s,\ q$ stand for nonnegative integers, and $g\geq 2$,
unless otherwise specified.

\subsection {Kodaira dimension and varieties of general type}
Given an abstract
algebraic variety
$V$, 
 consider its   rational images
 in projective space,  $\phi :V\dashrightarrow \P ^r$ where
$\phi$ is a rational map.
To such data one  canonically associates a line bundle $L$  on $V$: 
$L:=\phi ^* {\cal O}_{\P ^n}(1)$. Viceversa, given  a ``sufficiently
nice" line bundle $M$ on $V$ one constructs a rational mapping
$$
\phi _M:V\dashrightarrow \P (H^0(V,M)^*).
$$
Here are two definitions to clarify this set-up (see \cite{KM} for more details).

\begin{defi} Let $L$ be a line bundle on $V$
\begin{enumerate}
\item  $L$  is {\it big} if for some 
positive integer $k$ the map  $\phi _{L^k}$ gives an embedding 
 in projective space of a
non-empty open subset of $V$.
\item  $L$  is {\it ample} if for some 
positive integer $k$ the map  $\phi _{L^k}$ gives an embedding  of $V$  
 in projective space.
\end{enumerate}
\end{defi}
Assume now that $V$ is nonsingular and consider  its canonical line bundle, $K_V$. 
If, for a given
integer 
$n\geq 1$, the rational mapping 
$$\phi _{K_V^n}:V\dashrightarrow \P (H^0(V,K_V^n)^*),$$
is defined,   it will be  called
the
{\it $n$-canonical mapping} of $V$. 

An important  birational invariant
is the {\it Kodaira dimension} $\kappa (V)$
$$
\kappa (V) :=\max _{n\geq 1} \{ \dim {\rm Im} \phi _{K_V^n}\}
$$

It might very well happen that $K_V$ and
its powers have  no sections 
(for example, this will be the case for $V$ birational to projective space);
in such a case the Kodaira dimension of $V$ is set
to be equal to $- \infty$ (so that $\kappa$ is additive with
respect to products).

Of course we always have $\kappa (V)\leq \dim V$. We shall be particularly interested in the case
when  equality
holds, that is, 
$K_V $ is big.

\begin{defi}
\begin{enumerate}
\item A nonsingular variety $V$  is {\it of general type} if $\kappa (V)= \dim V$ 
\item A  variety  
is of general type if a desingularization  of it is of general type. 
\end{enumerate}
\end{defi}
 Part (2)
 is   well defined, since the Kodaira dimension is a birational invariant. 

On the
other hand, the  definition
presents  a subtlety that will play a crucial role in the sequel (see
section~\ref{sem}). 
Let $V$  be a (possibly singular) variety possessing the dualizing line bundle
$\omega _V$. Recall that $\omega _V$ is the unique line bundle making Serre duality valid
and, if $V$ is nonsingular, $\omega _V=K_V$. Now the warning  is:
 to say that a
singular variety
$V$ is of general type is not the same as saying that $\omega _V$ is big.
Consider in fact the following 
\begin{example}
 \label{quartic} $V\subset \P ^2$ is an irreducible
curve of degree $4$ having $3$ nodes. Then 
the dualizing sheaf of $V$ can be computed by the
adjunction formula,
$$
\omega _V = K_{\P ^2}\otimes {\O }_{\P ^2}(V)\otimes \O _V= \O _V (1)
$$
hence  $\omega _V$ is the line bundle defining the embedding of $V$ in $\P
^2$, in particular $\omega _V$ is big.
On the other hand,  
there is a birational morphism
$r :\P ^1 \la V$;  to obtain it, consider
the
$\P ^1$ of all conics $C_t$, with $t\in \P ^1$, passing through the $3$ nodes of $V$ and any other
fixed point of $V$, then by B\'ezout's theorem $C_t$ meets $V$ in one remaining point, thus determining
the morphism $r$. Therefore
$\P ^1$ is the desingularization of $V$, hence $\kappa (V)=-\infty$ and 
 $V$ is not of general
type.
The point  is that although
there are plenty of differential forms that are regular
on the complement of the three nodes of $V$ (which follows from $\omega _V$ being ample),
 the
pull-back of these forms via the morphism $r$  will not extend to regular forms on the
whole of  $\P ^1$. The singularities of $V$ are, in a suitable sense, too complicated (with
the appropriate terminology we shall say that they are not ``canonical", see section
\ref{sem}).
\end{example}

Some more useful examples. 

\begin{example}
 \label{gtcu} A curve is of general type if and only if its 
 genus is at
least $2$; if the curve is singular, by ``genus"  we mean the genus
of its normalization. 

\end{example}

\begin{example}\label{gtCI} Smooth complete intersections of high enough degree in
projective space are of general type; more precisely if $V\subset \P ^r$ is
the complete intersection of $c$ hypersurfaces $H_1, ....,H_c$ with $\deg
H_i=d_i$ (and $\dim V = r-c$) then $V$ is of general type if and only if
$\sum d_i \geq r+2$ (this fact follows directly from the adjunction
formula).  

\end{example}

\begin{example}\label{gtqu} If $V$ is of general type and $G$ is a subgroup of automorphisms of $V$
(so that $G$ is finite), then $V/G$ needs not be of general type,
but for $m$ large enough the quotient $V^m/G$ by the diagonal action of $G$
on $V^m$ is of general type. Notice that this is true regardless of $V$
being smooth (see \cite{CHM}     corollary 4.1).

\end{example}

Conjecturally, the property of being of general type should be strictly related to the
analytic, differential, and diophantine structure of a variety.

\subsection{Moduli of algebraic varieties}
\label{mod}

Fix a polynomial $h(x)\in \Q [x]$. 
Denote by $M_h$ the moduli
scheme of smooth algebraic varieties whose canonical bundle is ample and has Hilbert polynomial equal to
$h(x)$.  

The
simplest  case is when $h$ has degree $1$, 
then our moduli space  is non-empty if there is an integer $g\geq 2$ such that
$$
h(x)=(2g-2)x-g+1
$$
and  $M_h$  will be the moduli space of curves
of fixed genus
$g$, which is usually denoted by $M_g$.

For the main properties and the construction of such moduli spaces we refer to
the book of E. Viehweg \cite{Vibook}.

We need to recall here the fact that such spaces $M_h$ are not complete,
and in very few cases modular completions for them exist.
One of such cases is $M_g$, which is an integral, normal quasiprojective variety,
admitting a modular compactification, $\Mg$, which is the moduli space of Deligne-Mumford
stable curves; $\Mg$ is a
 normal, integral, projective variety.

If $V$ is a smooth algebraic variety as above, we denote by $[V]\in M_h$
the point of the moduli scheme corresponding to $V$.

Recall that if $f:X\la B$ is  family of nonsingular varieties  as above,
then there exists a unique regular morphism, called the {\it classifying morphism}
$$
\begin{array}{lcc r}
\gamma _f:& B &\la& M_h \cr
\;&b &\mapsto &[X_b]
\end{array}
$$

More generally,
if $X\la B$ is a family such that the locus $U\subset B$ where the fiber is smooth is open,
we shall continue to denote by $\gamma _f: B \da M_h$ the rational map
whose restriction to $U$ is the  classifying morphism.

\begin{defi}
\begin{enumerate}
\item
 A family $X\la B$   is  {\it isotrivial} if there exists 
a dense open subset $U\subset V$ over which all fibers are isomorphic.
\item A family has {\it maximal variation of moduli}
if there is a non-empty open subset $U\subset B$ such  for every $u_0\in U$ the set
$
\{u\in U : X_u\cong X_{u_0} \}
$
is finite.
\end{enumerate}
\end{defi}

In particular, a family of generically  smooth varieties for which there exists a moduli scheme
is isotrivial if and only if its classifying map
is constant, and it has maximal variation of moduli if and only if its classifying map is generically
finite. 

In fact, moduli theory allows us to measure how much a family varies in moduli.

\begin{defi}
Let $f:X\la B$ be a family of generically smooth varieties of general type.
Then the {\it variation of moduli of the family},
denoted $Var(f)$ is the dimension of the image of its classifying map,
that is, the number
$$
Var(f) := \dim \Im \gamma_f
$$

\end{defi}
Such a definition may appear to make sense only if there exists a moduli scheme
(e.g. if the canonical bundle is ample), but in fact
it works more generally, 
 as  stated (see section 2 of \cite{Kol}).

\section{From geometry to arithmetic}
\label{dio}
We shall assume throughout this section that $V$ is nonsingular.

\subsection{The one dimensional case}
\label{one}
For curves 
the  picture is  
very suggestive: the algebro geometric condition of being of general type
(that is, of having genus at least 2) can be characterized in terms of

(1) topology

or (2) complex analysis

or (3) differential geometry

or (4) arithmetic.

\ni
More precisely, here are some  well known fundamental facts. 
 
\

\ni
{\bf 1.} The universal covering of $V$ is the open unit disc
(Uniformization).

\ni
{\bf 2.} If $f:\C \la V$ is a holomorphic map, then $f$ is constant (Picard-Borel
theorem).

\ni
{\bf 3.}  $V$ carries a metric of constant negative curvature.

\ni
{\bf 4.} $V(K)$ is finite for any number field $K$ (Mordell conjecture - Faltings theorem).

\

\ni
Conversely, each of the properties listed
above  implies that the curve is of general type.
In fact, if $V$ has genus $0$ then $V$ is the Riemann sphere.
Hence $V$ 
is
simply connected, it contains copies of $\C$ as algebraic open subsets and it   carries a
metric of constant positive curvature. Finally it is  easy to see that $V$ has infinitely
many rational points over a sufficiently large number field.

If $V$ has genus $1$, then $V=\C /\Gamma $ is a complex torus and the
quotient map 
$\C \la \C /\Gamma$ is the universal covering
(with the intrinsic identifications $\C \cong H^0(V,K_V)^*$
and $\Gamma \cong H_1(V,\Z)$); moreover  $V$ carries a
flat metric. If $V$ contains a $K$ rational
point, then $V(K)$ can be given a group structure. The theorem of Mordell-Weil states that
the set of $K$-rational points is a finitely generated abelian group 
$V(K)=\Z ^r\oplus H $ with $H$ a finite group. 
It is well known that for any curve of genus 1 (or for any abelian variety)
there exists a  finite extension $F$  of $K$, such that $V(F)$ is
infinite.

Additionally, in the genus 1 case, natural  uniformity problems  arise, such as: 
(1) How can $H$ vary? (2) Is  $r$ uniformly bounded, or for any
$N$   there is an elliptic curve  for which $r>N$?

Question (1) has a definite answer thanks to a theorem of Mazur
('76) (for $K=\Q
$) showing that there exist exactly sixteen
 finite groups (including the trivial group) that do occur as
the torsion part of $V(K)$. 
This has been   generalized over all number fields by
Merel \cite{Me}  who proved that the number of torsion points of an elliptic
curve defined over a degree $d$ extension of $\Q$ is uniformly bounded.
The second question is  open.

\subsection{Higher dimensional varieties} 
Thus, in the world of curves 
there is a  general dichotomy between curves that are of general type and curves that are
not. It is natural to ask if something similar
 occurs for varieties of higher dimension. So,  let
$V$ be a
variety of general type; what can one
 say about $V$ from a complex analytic, differential geometric,
diophantine point of view? 
With the one-dimensional case in mind,  recall the  analytic terminology:

\begin{defi}
A variety $V$ is 
\begin{enumerate}
\item
{\it Brody hyperbolic}
 if  any holomorphic map $f:\C \la V$ is
constant. 
\item
{\it pseudo-hyperbolic} if the locus of all images of
holomorphic, non-constant maps $f:\C \la V$ is not Zariski dense.
\end{enumerate}
\end{defi}

It is clear that  Brody hyperbolicity is not  birationally invariant,
which motivates the second 
part  of the above definition. 

Turning now to differential geometry: given the hyperbolic distance
 on the open unit
disc one  defines the Kobayashi semi-distance $d_{kob}$ on $V$
(see \cite{L} for details).
Then 
\begin{defi} 
A variety 
$V$ is  {\it Kobayashi hyperbolic} when $d_{kob}(x,y)=0$ 
implies $x=y$, that is, if $d_{kob} $ is a distance on $V$.
\end{defi}
It turns out that a holomorphic map is always distance decreasing with respect to
the Kobayashi semi-distance, moreover, 
the Kobayashi semi-distance is identically $0$ on $\C$. Hence if
$V$ is Kobayashi hyperbolic, then $V$ is Brody hyperbolic. 

It is an
important theorem of Brody that 
the converse is  true (provided that $V$ is compact).
Hence  the  analytic   and the differential geometric notions of hyperbolicity
coincide
for projective varieties.

The expected link with algebraic geometry is 
\begin{conj}[Lang \cite{L}] $V$ is hyperbolic if and only if  $V$ and all of its
subvarieties are of general type
\end{conj}

We briefly describe two non trivial examples  of varieties of general type which are
not  hyperbolic.
The first (of Mike Artin) is a 
smooth surface of degree
$5$ in $\P ^3$ (whose canonical bundle is  very ample, by  \ref{gtCI}) 
having a plane
section $C$ with $6$ nodes, so that $C$ is a rational curve.
The second  is a whole class of examples: the Fermat's 
hypersurfaces of degree at least $r+2$ in $\P ^r$. They always contain lines,
and, again by  \ref{gtCI},  they are of general type.

Kobayashi suggested that such examples should be exceptional, 
by conjecturing (in \cite{Kob}) that a generic
hypersurface of high enough degree in $\P ^{r+1}$ is hyperbolic; 
it is known that, for such a conjecture to hold in dimension $r\geq 2$, the degree
would have to be at least $2r+1$.
In '98 one  important step was taken by Demailly-El Goul \cite{DE}, 
who proved the
 conjecture  true for very generic surfaces in
$\P ^3$ of degree at least
$42$; here the terms ``generic" and ``very generic" mean that the exceptional locus  is a
finite, respectively countable, union of proper closed algebraic 
subsets of the moduli space of hypersurfaces.

The significance of such a Theorem is better understood 
when compared with a previous famous result
of Clemens and Xu, stating that a very generic surface in $\P ^3$ of degree at least $5$
contains no algebraic curves of genus less than $2$, leaving the case of
transcendental curves open.

The birational  analogue is the following well known

\begin{conj}[Green-Griffiths \cite{GG}]
If $V$ is of general type then $V$ is pseudo
hyperbolic.
\end{conj}

If the above statement were true, it would imply that, in a variety of general
type, the union of all
curves of genus less than $2$ is not Zariski-dense. In particular, a surface of general type would contain  only finitely
many curves of genus $0$ and $1$. This is  known to be true 
for  surfaces satisfying the  inequality
$c_1^2>c_2$; Bogomolov proved this for algebraic curves (\cite{B}),  his method was
 extended to holomorphic curves by McQuillan
(\cite{Mc}). Some more progress about the above conjecture in dimension
$2$ is in
\cite{LM1}.

\subsection{Diophantine conjectures.}
\label{ffdio}
What about arithmetic properties of varieties of
general type? 

The main conjectures on the distribution
of rational points  in varieties of general type have been partly
motivated by parallel results over function
fields. 

To start with, the  analogue of
the Mordell conjecture  for function fields
was proved by Manin \cite{Man} about twenty years before the proof over number
fields (by Faltings). More precisely, let  $L$ be the function field of a 
projective variety $B$ over the complex
numbers.
Manin showed that a curve defined over $L$, having genus at least $2$,
either has finitely many $L$-rational points, or it is birationally
a product and all but finitely many of its $L$-points correspond to
constant sections.

Besides Manin's, there have been other approaches to the geometric Mordell problem.
We are mostly interested in one of them   which on the one hand
 highlights the analogies with the number field case; on the other hand establishes
a link with moduli theory.
We begin with some notation

\begin{defi}
\label{curves}  Let $F$ be a number field or a   function field $F=\C(B)$.
\begin{enumerate}  \item Denote by $\CF$ the set of $F$-isomorphism classes of
nonisotrivial,  nonsingular, projective  curves of
genus
$g$ defined over $F$. 
\item 
Let $S$ be either (if $F$ is a number field) a finite set of places of  $F$ or, if $F=\C(B)$,
a closed subset of $B$. 
Denote by $\CFS$ the subset of $\CF$ parametrizing curves having good reduction outside of $S$
\end{enumerate} 
\end{defi}
(In the aritmetic case the   nonisotriviality condition is 
of course useless.) 

During the ICM in  1962, Shafarevich  conjectured that
{\it If   $F$ is
either a number field or the function field of a complex  curve,
then $\CFS$  is finite. }

The geometric case  was obtained by
Parshin (if $S$ is empty) and Arakelov (in general) ( see \cite{P}, \cite{Ar});  they proved that
 if $B$ is a smooth complex curve and $S\subset
V$ a finite subset, then there exist only finitely many non-isotrivial families of smooth
curves of genus $g\geq 2$ over $B-S$. 

The connection with the Mordell problem was at the same time established by Parshin, 
in \cite{P}, where, by the celebrated ``Kodaira-Parshin trick",  he showed that the
Shafarevich conjecture implies the Mordell conjecture, both for function fields and
for number fields:
\begin{lemma}[Kodaira-Parshin trick]  Let $F$ be a number field or the function field of a complex curve.
Assume that  the set $\CFS$ is finite for every $F$ and  $S$ as above.
 Then  for every $C\in \CF$ the set $C(F)$ is
finite
\label{trick}
\end{lemma}

The number field version of the Shafarevich conjecture (and hence of the Mordell conjecture)
was proved
by Faltings in \cite{F}.

Going back to the Theorem of Manin,  it turns out (\cite{Man} proposition 5) that
the case of higher dimensional bases $B$ can be reduced, by induction on the dimension,
to the case $\dim B =1$, to which  we therefore restrict ourselves.
Then the
second part of the theorem says that, if the family is birational to a product
$B\times C_0$, where $C_0$ is a curve of genus at least $2$ defined over
$\C$, then there can only be finitely many non-constant sections $B\to
B\times C_0$. This  is a straightforward consequence of a classical
theorem of De Franchis (\cite{DF}), according to which there can only be finitely many
non-constant morphisms from $B$ to
$C_0$, provided that $C_0$ has genus $2$ or more.

In '75 Kobayashi and Ochiai generalized De Franchis Theorem to
complex varieties of any dimension:

\begin{thm}[\cite{KO}]
\label{KO} Let $Y$ and $Z$ be projective  varieties of the same dimension, with
 $Y$  of general type. Then
there exist  finitely many dominant rational maps from  $Z$ to
$Y$. 
\end{thm}
A generalization of Manin's theorem   to varieties of higher
dimension was obtained by Noguchi\cite{N}:
\begin{thm}[\cite{N}]
Let $X$ be smooth of general type defined over a function field $L$.
Assume that $X$ has ample cotangent bundle. Then either $X(L)$ is not
Zariski-dense, or $X$ is definable over $\C$ and all
but finitely many of its $L$-rational points are defined over $\C$.
\end{thm}

Clearly, if $X$ has ample
cotangent bundle, then $X$ is of general type and so are all of its
subvarieties. The converse is false.

Such a set  of  results and trends  
inspired the following speculations in the arithmetic setting.
Let $K$ be a number field and $V$ a variety defined over $K$.

\begin{conj}[Weak form]
If $V$ is of general type, then $V(K)$ is not Zariski-dense in $V$.
\end{conj}

This statement is  attributed to E. Bombieri,  S. Lang and P. Vojta, and often referred to as
the ``Weak Lang conjecture". The next stronger form is due to Lang.

\begin{conj}[Strong form]
If $V$ is of general type, then there is a proper, closed subset
$V_0$ of $V$ such that, for any number field $K$, all but finitely many
$K$-rational points of $V$ lie in $V_0$.
\end{conj}

The above conjectures are  known to be true only for curves (\cite{F}) and subvarieties
of abelian varieties (\cite{F2}). 

We will talk more about them (namely, about their implications) in the sequel.

\subsection{When are rational points Zariski dense?}  Consider  varieties
of low Kodaira dimension and look for geometric conditions that should ensure density of
rational points over some number field. 

We start with some terminology: 

\begin{defi} Let $V$ be a variety defined over a number field $K$; rational points of $V$ are {\it potentially
dense}  if there exists a finite extension
$K'$ of $K$ such that
$V(K')$ is Zariski dense in $V$.
\end{defi} 

According to the diophantine conjectures, rational points should not be  potentially dense on varieties
admitting a dominant map to a variety of general type. What about the rest of the world?
A result of Colliot-Th\'el\`ene, Skorobogatov, Swinnerton-Dyer \cite{CSS} shows that it is
wrong to conjecture that rational points are potentially dense everywhere else.
They prove (Section 2 in \cite{CSS}) the following:

\begin{thm}\cite{CSS} Let $p:X\la \P ^1$ be a fibration with $X$ smooth, irreducible
and projective and such that the generic
fiber of $p$ is smooth. Assume that all objects above are defined over the number field $K$. 
If $p$ has at least $5$ double (i.e. non reduced) fibers, then the $K$-rational points of $X$
must all lie in a finite union of fibers.
\end{thm}
In particular, rational points  are not potentially dense on $X$. Explicit examples of such
fibrations are found already in dimension $2$, with $X$ a hyperelliptic surface which does not admit a
dominant map to a variety of general type (Proposition 3.1 in \cite{CSS}).

When are then rational points potentially dense?

Obvious or well known cases where potential density holds are rational, unirational and abelian
varieties. As ventured by J. Harris and Y. Tschinkel in \cite{HT}, natural  candidates are
varieties whose canonical bundle is ``negative".
They propose
\begin{conj}[\cite{HT}] Rational points are potentially dense on
\begin{enumerate}
\item all Fano varieties.
\item (stronger form) all varieties $V$ such that $-K_V$ is nef.
\end{enumerate}
\end{conj}
where  a smooth variety $V$ is called  
``Fano"   if
$-K_V$ is ample, and a line bundle $L$ is called ``nef" if for every curve $C$ in $V$,
$\deg _CL\geq 0$.

As a first piece of evidence, they
show that rational points are potentially dense on smooth quartic hypersurfaces of
dimension at least
$3$. The case of dimension $2$, quartic surfaces in $\P ^3$, is not
yet completely understood: in \cite{HT}  potential density is proved only under the additional 
assumption that the quartic surface contains a line
$\ell$. If this
happens the surface has  a fibration in elliptic curves, given by  all planes containing $\ell$, 
whose
intersection with the surface is given by $\ell$ union a varying residual plane cubic;  this is
what allows them to prove density of rational points. 

This example is quite typical: in almost all cases in which potential density is known, the method of
proof involves exhibiting a dense collection of subvarieties having a dense subset of rational points
over the same number field.

This is done in \cite{BT1} and \cite{BT2} for elliptic K3 surfaces and for Enriques surfaces
by means of elliptic fibrations, and in \cite{Si} for K3's having infinitely many
automorphisms. 

Part $2$ of the above conjecture is a special case  of Conjecture $III_A$
in \cite{Camp}, which proposes that {\it rational points are potentially dense on $V$
if and only if $V$ is  special} (in the sense of F. Campana: Definition 2.1 in
\cite{Camp}). 
The relation comes from the prediction  that varieties with nef anticanonical
bundle ought to be special (see \cite{Camp}, 2.2).

To conclude: not surprisingly, the results in dimension 2 and 3 owe a good deal
to the corresponding classification theorems; we summarize them  as follows.

Dimension $2$: 
potential density is known to hold for all surfaces which are not of general type,
with the exception of K3's with finite automorphism group and without elliptic fibrations,
which case is not known.  However, it is worth noticing that for every K3 surface, potential
density has been proved for its symmetric products of suitable order (\cite{HasT}).

In dimension $3$,  potential density is known for all Fano threefolds, with the
exception of double covers of $\P ^3$ ramified in a smooth sextic surface (\cite{HT}, \cite{BT3}).
It is unknown for the remarkable case of Calabi-Yau varieties.

\section{\bf  Moduli theory and uniformity conjectures}
\label{AG}

\subsection{Consequences of the diophantine conjectures}
If the diophantine  conjectures (stated in \ref{ffdio}) were true,
they would have  a big impact 
 about
the  distribution of rational points on
curves. 
To explain how, let us begin by applying them to a special case.
Let $f:X\to B$ be a family of curves of genus at least $2$, assume that 
$B$ is a nonsingular curve and that everything is defined over $K$.
 Faltings Theorem ensures that for every point $b\in B$, the fiber
$X_b$ of $f$ over $b$ has a finite
 set of  $K$-rational
points; consider the union $F_K$ of all such sets: 
$$F_K := \cup _{b\in B(K)} X_b(K)$$
Now ask whether the restriction of $f$ to $F_K$ has bounded fibers over $B$;
notice that $F_K$ might very well be dense in $X$.

In the special case in which  $X$ is 
 of general type,
the Weak Diophantine conjecture implies that there is a proper closed subset $Z$  of $X$ containing
$X(K)$; clearly then
$F_K\subset Z$. Thus, the restriction $f_{|Z}:Z\la B$ having finite degree
(since $\dim X =2$ and hence $\dim Z =1$), the fibers of $F_K$ over $B$
have bounded cardinality.

Of course, assuming that $X$ be of general type was  crucial. 
Now, let us weaken such a ridiculously restrictive assumption
by  the following hypothesis: given a family $f:X\la B$ of curves of genus at least $2$,
there exists a number $n$
such that the $n$-th fibered power 
$X\times _B  ....\times _B X$  has a dominant  map to a variety of general
type.
 Then a  subtle elaboration 
(\cite{CHM}, proof of Theorem 1.1)
of the quick reasoning above shows that, if this hypothesis
holds,  then the Weak Diophantine conjecture implies that the
sets of rational points of the fibers of $f$
have bounded  cardinality.

At this point the question becomes: when is the above  ``fibered power hypothesis" satisfied?

The answer is: always! That is, for every family of generically smooth curves of genus at least $2$.

The proof of this
algebro-geometric fact,  called the ``Correlation Theorem",
is far from being trivial and  will be dealt with in the next section.

Here we shall focus on its arithmetic consequences. 
They are obtained by applying the above machine
to a special  family of curves: a so-called ``global family",
that is a morphism $f:X\la B$ such that every curve of genus $g$ is isomorphic to some
fiber of $f$. It is well known from moduli theory
that there are many choices for such a global family, and that one can assume
$B$ to be nonsingular.

The first corollary of the Correlation Theorem is   that  if the Weak Diophantine
conjecture  holds, then there exists a uniform bound $\ub <\infty$ such that
any curve of genus $g\ge 2$ defined over $K$ has at most $\ub$
$K$-rational points (\cite{CHM}).
This result has been  strenghtened  by Abramovich and Pacelli,
who showed that such a bound  only depends on the degree $d$ of $K$
over $\Q$. 
\begin{thm}
\label{ub}
The  Weak Diophantine Conjecture
implies that for any $d\ge 1$ and for any $g\ge 2$ there is a number
$B_g(d)$ such that for every $K$ with $[K:\Q ]\le d$ and for every curve $V$ of
genus
$g$, defined over $K$ we have
$|V(K)|\le B_g(d)$.
\end{thm}
(This is Theorem 1.1 in \cite{CHM} with strengthening from \cite{A1} and \cite{Pa}.)

In a  similar vein the Strong Diophantine conjecture implies
a  universal uniformity statement (Theorem 1.2 in \cite{CHM}):

\begin{thm}
\label{uub}
The  Strong Diophantine Conjecture
implies that for any $g\ge 2$ there is a number $\gb$
 such that for every $K$ there exist only finitely many curves having more
than $\gb$  points defined over $K$.
\end{thm}

Of course  the (finite)
set of curves having more that
$\gb$ $K$-rational points has to depend on $K$. 

At present, it is not known whether either of the two above statements, following from the diophantine 
conjectures,
hold.  A true consequence  was obtained by Mazur (in \cite{Maz}):
\begin{cor}
\label{iso}
The Strong Diophantine Conjecture implies that for any $K$ there are only finitely
many $j$-invariants of elliptic curves defined over $K$ and having more
than $20$ $K$-rational isogenies.
\end{cor}
For $K=\Q$ the above  finiteness result  is  known to be true independently of the Diophantine
Conjecture:  a theorem of Kenku \cite{Ke} states that any elliptic curve over $\Q$
has at most $8$ rational isogenies.

\

\subsection{The role of algebraic geometry}
\label{cor}
As already mentioned, the proofs of the two Theorems  \ref{ub} and \ref{uub}  rely on an
algebro-geometric result,   the Correlation
theorem, which says that  high fibered powers of a family of curves of genus at least
$2$ dominate a variety of general type
(see below).
The name  Correlation
theorem was given to it  because of the application involving the Diophantine
conjecture:  if $X\to B$ is a family of curves,
then by the Correlation theorem  there exists a variety of general
type $V$ such that the $n$-th fiber power
$X\times _B .... \times _BX$ has a (rational) dominant map $h$
to $V$. Now
the Weak Diophantine conjecture says that the rational points of $V$ are contained
in a proper, closed subset $V'$ of $V$;
the preimage 
of $V'$ via $h$ will be a closed subset $Z$ of the fibered product, 
containing all of its  rational points. So the algebraic functions defining
$Z$ can be viewed as co-relating functions for $n$-tuples of rational
points of the fibers of $X\to B$.

In \cite{CHM} the Correlation theorem was proved   for 
curves, using their moduli theory;
it was  conjectured to hold for families of varieties of any dimension. 
From a purely algebro-geometric point of view, the issue was 
 interesting in its own right (regardless of its applications to arithmetic), and fit
very well with a lively area of research, at the time.
In fact it did not take long for its full generalization  to
be obtained by D. Abramovich.
Here it is,  with a more appropriate name (\cite{A3}):

\begin{thm}[Fibered Power Theorem]
\label{fpt}
Let $f:X\to B$ be a proper morphism of integral  varieties such that the general
fiber is  smooth  of general type. Then there exists an integer $n$ such
that
$X^n_B$  has a dominant rational map to a variety of general type.
\end{thm}

Where  $\Xn$ is the irreducible component of the
$n$-th fiber product of $X$ over $B$ 
which dominates $B$; we shall denote
$$f_n:\Xn \la B$$
the natural map (so that $X = X^1_B$ and $f=f_1$).

\begin{remark} The proof of the theorem   explicitely constructs a variety of general type dominated by $\Xn$.
Call it $V$, then  if $f:X\la B$ is defined over a number field $K$,  so is $V$. In fact
$V$ is tightly related to  the given family, for example, if $f$ has maximal variation of moduli
(see section~\ref{mod}), then $V$ is equal to $\Xn$
(theorem~\ref{maxvar}). In general, if the fibers of $f$ have dimension $N$, then 
$$
\dim V = Var(f) + n + N
$$
($Var(f)$ being the variation of moduli of $f$, defined in \ref{mod})
\end{remark}

 B. Hassett, in
\cite{Ha}, proved
it for families of surfaces, by applying the moduli theory
of surfaces of general type; Abramovich obtained the general
case
 combining an argument
suggested by Viehweg (using results of Koll\'ar) together with some
techniques of de Jong. The cases of curves and surfaces were 
dealt with by making use of semistable reduction and
of existing compactifications for
their moduli spaces.
Such tools were not available for higher dimensional varieties.

We give now an outline of the proof of this theorem, comparing the original
argument of \cite{CHM} with that in \cite{A3}.
Without further mention we will assume that all maps
are  rational maps, unless we specify otherwise; we will use the notation
$X\bir Y$ to mean $X$ birational to $Y$.

\begin{defi}
A morphism
 $f:X\la B$ is a {\it family of relative general type}
if $f$ is a proper  morphism of integral
varieties such that the generic fiber of $f$ is  smooth and of general type.
\end{defi}

The key result  is the special case of ``truly varying" families (see section \ref{mod})

\begin{thm}
\label{maxvar}
Suppose that $X\la B$ is a family of relative general type
having maximal variation of moduli. Then for $n$ large, $\Xn$ is of general type.
\end{thm}

Then  the Fibered Power Theorem is proved by  reducing it to the case of a
family with maximal variation of moduli. This is done by means of
  the following diagram, whose  existence follows from standard
moduli theory

\

$$
\begin{array}{lcccccccr}
X_1/G & \bir & X & \longleftarrow & X_1 & \bir & T\times _S B_1 & \la & T\\
   \; & \; & \downarrow & \; & \downarrow & \; & \downarrow & \; & \downarrow \\
   \; & \; & B & \longleftarrow & B_1  & = & B_1 & \longrightarrow & S  \\
\end{array}
$$

\

Where  the starting point
is a family $X\la B$ of relative general type, as in the statement of Theorem~\ref{fpt};
$B_1\la B$ is a generically finite surjective Galois morphism with group $G$,
and
$$
X_1:=X\times _B B_1
$$
$G$ acts as a subgroup of $\Aut _{B_1}(X_1)$ in a way that $X_1 /G $ is birational to $X$.
The map $B_1\la S$ is a   classifying map
(identifying points in $B_1$ with isomorphic fibers) and $T\la S$ a family with maximal variation of
moduli, which should be thought of as  a tautological, or universal, family
(in particular, for every point $s\in S$ the fiber of $T$ over $s$ is isomorphic to
the fiber of $X_1$ over any point in $B_1$ that maps to $s$).
Finally,
$X_1$ is birational to $T\times _S B_1$ over $B_1$.

The diagram is such that for every $n\geq 1$ there is a rational dominant map $\Xn \la T^n_S/G$
($G$ acting diagonally). By Theorem~\ref{maxvar}, $T^n_S$ is of general type for large $n$,
therefore, by \ref{gtqu}, 
$T^m_S/G$ is also of general type for $m$ large,  hence the
fibered power theorem is proved.

The simplest nontrivial example of how this  works, is that of an
isotrivial family: assume that for every $b$ in some open subset of $B$ we have that
$X_b\cong X_0$, where $X_0$ is a fixed nonsingular variety of general type.
Then, since $X_0$ has finitely many automorphisms (by Theorem~\ref{KO}) one can construct 
the finite covering $B_1\la B$ 
so
that the base changed family is trivial:
$X_1\bir X_0\times B_1$. The space $S$ is just one point and $T\cong X_0$ so that 
$X$ is birational to $(X_0\times B_1)/G$; $G$ is a subgroup of $\Aut F$ such that 
$X$ dominates $X_0/G$ and therefore, for every $n$,
$\Xn $ dominates $X_0^n /G$.

\subsection{Semistable reduction and its variants}
\label{sem}
We are then left with the proof of theorem~\ref{maxvar}.
The information that the fibers are of general type is contained in the relative
dualizing sheaf  $\omega _f$ of the family, which is relatively big, that is, 
$\omega _f \otimes {\cal O}_{X_b} = \omega _{X_b} $
  is big for general $b\in B$. The  question is: how does $\omega _f$
behave as a line bundle on the total space $X$?
The answer is provided by the

\begin{prop}
\label{big}
Let $f:X\la B$ be a family of relative general type, with $X$ and $B$  smooth.
Assume that the family has maximal variation of moduli. Then 
\begin{enumerate}
\item $\omega _f $ is big.
\item For $n$ large, $\omega _{\Xn}$ is big.
\end{enumerate}
\end{prop}
\begin{proof} The first part is the difficult one; an ad hoc proof for
families of curves is outlined in   \cite{CHM}
(Lemma 3.2). The general case  follows  from deep results of Koll\'ar and Viehweg (\cite{Vi1} and
\cite{Kol}). 

The second part  is not hard, and follows from the first: notice to start that
$\omega _{f_n}$ is also big. To prove that $\omega
_{\Xn}$ is big; we  have the relation
$$
\omega _{\Xn}=\omega _{f_n} \otimes f_n^*\omega _B =p_1^*\omega _f \otimes
....
 p_n^*\omega _f \otimes f_n^*\omega _B
$$
where $p_i :\Xn \la X$ is the $i$-th factor projection.
Therefore, what could prevent $
\omega _{\Xn}$ from being big is  $\omega _B$.
On the other hand the above formula suggests that such an obstruction does
not depend on $n$. In fact  one proves that, for $n$ large enough,
the positivity of $\omega _f$ overcomes the negativity of $\omega _B$ so that
 $\omega _{\Xn}$  is big.
(Lemma 3.1 in \cite{CHM}).
\end{proof}
Our goal is to prove that $\Xn$ is of general type and we must now use caution,
because $\Xn$ may be singular:
as we saw in  \ref{quartic}, the fact that $\omega
_V$ is big for a singular variety $V$ does not imply that $V$ is of general type. 

 Let $r:V^r\la V$ be a resolution of singularities,  we
shall say that
$V$ has {\it canonical} singularities if $\omega _V^m$ is a line bundle for some integer $m$,
and if,  
given a regular section $\sigma$ of $\omega
_V^m$, the pull-back
$r^*\sigma$
 extends to a regular section of $K_{V^r}^m$ (see \cite{KM} 2.11).  

In particular, if
$V$ has canonical singularities and
$\omega _V$ is big, then 
$K _{V^r}$ is also big and $V$ is of general type.

What is needed now is an analysis of
the singularities of $\Xn$: if they were canonical, the second part of
the above Proposition would, of course, imply
that $\Xn$ is of general type.

We are
here facing the most delicate part of the proof,
 where  the argument in \cite{A3}  must be substantially different from that in \cite{CHM},
as we will explain.

The problem is that  we  have no control over the singularities of $\Xn$.
The solution is to relate the given  $X\la B$ to a new family 
whose  fibered powers have canonical
singularities. The crucial tool that will be used is:

\begin{prop}
\label{can}
Let $f: X\la B$ be a proper morphism of generically smooth integral varieties
such that locus in $B$ where the fiber is singular is a divisor with normal
crossings. Assume that one of the two following conditions hold
 \begin{enumerate}  
\item All fibers of $f$ have normal crossings singularities.
\item The family $X\la B$ is plurinodal
\end{enumerate} 
Then for every $n$ the singularities of $\Xn$ are canonical.
\end{prop}

Where 
a  family $f: X\la B$ is called plurinodal if 
$f$ factors through finitely many morhisms
$$
f:X= X_h\la X_{h-1}\la .... \la X_1 \la B = X_0
$$
where each $X_i\la X_{i-1}$   is a family of irreducible, nodal curves.

The proof of part 1  for curves is  Lemma 3.3 in \cite{CHM},
for the general case and part 2  see Lemma 3.6 in \cite{Vi1} and section 4.2 in \cite{A3}.

\begin{remark}
The normal crossings assumption on the degenerate locus of $B$ is always attainable after
birational modifications, which will not alter our final conclusions.
\end{remark}
The first   condition is used  if semistable reduction  is known, for example, for curves. More
precisely, it is well known that,
if $f:X\la B$ is a family of generically smooth curves,  there is a
diagram:

\

$$
\begin{matrix}
Z/G  \bir & X & \longleftarrow & X_1 & \bir &Z\\
 \;   &\mapdown{f} & \; & \mapdown{f_1}& \; & \mapdown{} \cr
 \; &B & \longleftarrow & B_1  & = &B_1 \cr 

\end{matrix}
$$

\

\ni
where $X_1=X\times _B 
B_1$. The new object  $Z\la B_1$ is the so called ``semistable reduction"
of $X\la B$, that is, a family of curves with normal
crossings singularities (that is, ordinary double points) which coincides with   $f_1: X_1\la B_1$ away from the
singular fibers of $f_1$. Moreover $B_1\la B$ is a generically finite, surjective Galois covering with group
$G$ and $G\subset \Aut
_{B_1} (Z)$ so that $X$ is birational to $Z/G$. 
Thus, if $X\la B$ has maximal variation of moduli, so does $Z\la B_1$.
By Proposition~\ref{can}
the fibered powers $Z^n_{B_1}$ have canonical singularities; therefore, by Proposition~\ref{big},
 for large
$n$, 
$Z^n_{B_1}$ is of general type. Finally we get that,
for large $m$,  $X^m_B$ is of general type, being birational to the
quotient
$Z^m_{B_1}/G$ (using \ref{gtqu}).

The question remains as how to treat the higher dimensional case,
where  a  diagram like the above
does not  exist.  Abramovich's method uses techniques developed
in
\cite{De} and 
 replaces semistable reduction  by a  sort of plurinodal reduction.
Here is the diagram that he constructs.  (Lemma 4.1 in \cite{A3}):

\

$$
\begin{matrix}
X & \longleftarrow & X_1 \bir Y/G & \longleftarrow &Y\\
\mapdown{f} & \; & \mapdown{}& \; & \mapdown{h} \cr
  B & \longleftarrow & B_1  & = &B_1 \cr 

\end{matrix}
$$

\

\ni
where the new  family $Y\la B_1$  is 
plurinodal.
Furthermore, $B_1\la B$ is a generically finite, surjective morphism,
$X_1=X\times _B 
B_1$ and 
$G\subset \Aut _{B_1}(Y)$ is a finite group of automorphisms of $Y$ such
that $Y\la X_1$
 is the quotient map (that is, $X_1 $ is birational to $Y/G$).

Now,  condition 2 of Proposition~\ref{can} is obviously satisfied  by $Y\la B_1$, so that
for every
$n$ the singularities of
$Y^n_{B_1}$ are canonical. 

The trouble is that
$Y$ and the original family $X\la B$ are not as naturally  related as  in the  semistable reduction
diagram. The only connection is that  $Y/G$ 
(and similarly $Y^n_{B_1}/G$ for $G$ acting diagonally) is birational to the  base
change of
$X\la B$ 
(or of $\Xn \la B$)
by a generically finite map.

This last diagram is in fact harder to handle than the previous one: 
there are  a few  technical steps, which we omit, needed to descend the positivity
properties of
$\omega _h$ and $\omega _Y$ down to $X$ and show that, for $n$ large, $\Xn$ is of general type.

\subsection{Integral points.} For integral points on affine varieties, an argument similar
to the proof of theorem~\ref{ub}
yields that the analog of the Diophantine  conjectures implies certain
uniformity phenomena.

In 1929 Siegel proved that if $U$ is an affine curve whose projective closure has genus at
least $1$, then $U$ has finitely many integral points over any number field.

Such a theorem of Siegel has a conjectural generalization, due to Lang
and Vojta.

\begin{conj}[Lang-Vojta]
Let $X$ be a variety of logarithmic general type and let $\cal X$ be a
model for $X$ over $\Spec \O_{K,S}$. Then the set of $S$-integral points,
$\cal X(\O_{K,S})$, is not Zariski-dense.
\end{conj}

Where $S$ is a finite set of places of $K$ and

\begin{defi}
A quasiprojective variety $U$ is of {\it logarithmic general type} if given
a smooth compactification $V$
of a desingularization $U'$ of $U$, such that $D:=V\setminus U'$ is a a
divisor with normal crossings, then, for $n$ large, $n(K_V + D)$ gives a
birational map of $V$ in projective space.
\end{defi}

For our purposes it suffices to know that 
\begin{enumerate}
\item A projective variety is of logarithmic general type
if and only if it is of general type 
\item A curve of positive
genus with one or more points removed is of logarithmic general type
\item  $\P ^1$ with $n$ points removed is of logarithmic general type iff
$n\ge 3$. 
\end{enumerate}
Abramovich extended the techniques of \cite{CHM} to the case of integral points
on elliptic curves. He proved that the Lang-Vojta conjecture implies
uniformity of integral points on certain special models of elliptic curves.
One of his results is the following (see \cite{A2} Corollary 1):

\begin{thm} 
The Lang-Vojta conjecture implies that the number of integral points on
semistable elliptic curves defined over $\Q$ is uniformly bounded.
\end{thm} 
In fact, he proves a more general result, for any number field $K$
and any
set of places $S$.

Finally, just as in the case of rational points, Pacelli (\cite{Pa1})
showed that there exist a  uniform bound for $S$
integral points on semistable
elliptic curves which only depends on the degree of $K$ over $\Q$ and on $S$.

The uniformity consequences of the diophantine conjectures that we listed in this survey are 
open,
with the exception of corollary~\ref{iso} and of a remark of  Abramovich (unpublished)
who showed  that a true consequence of
Lang-Vojta conjecture is the already mentioned 
uniform boundedness of torsion points on elliptic curves,
  proved in \cite{Me} (see \ref{one}).

\

\section{Boundedness results}

\subsection{Counting rational points on curves.}
\label{count}
This section is devoted
to     uniformity questions about rational points on curves,
forgetting from now on  the Diophantine conjectures.

We  ask: what are the best known bounds on $\ub$ and $\gb$ (see Theorems~\ref{ub} and \ref{uub} for the
definitions of
$\ub$ and $\gb$)? 
In other words: 

\begin{enumerate}
\item (for $\ub$) How many rational points can a curve have? 

\item (for $\gb$) How many rational points can an infinite collection of curves
(all of  genus $g$ and defined over the same number field) have? 
\end{enumerate}

It would obviously be optimal  to have upper bounds for them, but this
 stands completely open as of today; in particular, $\ub$ and $\gb$ might
very well be infinite (which would of course imply that the diophantine conjectures are
false!). 

We then investigate lower bounds, that is, we look for examples of
curves having as many rational points as possible.
 Until '93, before the results described earlier revived interest in this
subject, the records for small genus, over $\Q$,
were held by Brumer who had  examples  showing that
$B_2(\Q )\geq 144$ and $B_3(\Q )\geq 72$. Since then, these records have
been broken, and as of today one knows that
$$
B_2(\Q )\geq 588
$$
with the following curve discovered by Kulesz:
$$
y^2 = 278271081x^2(x^2-9)-229833600(x^2-1)^2
$$
which has $12$ automorphisms.
For genus $3$,
$$
B_3(\Q )\geq 112
$$
with the hyperelliptic curve
$$
y^2=48397950000(x^2+1)^4-939127350499(x^3-x)^2
$$
found by Keller and Kulesz. This last curve has at least $16$
automorphisms; there are also examples of curves having many rational
points, but no extra-automorphisms; for instance the curve below was found by Stahlke,
and has $306$ points over $\Q$
$$
y^2=9703225x^6-9394700x^5+152200x^4+1124745x^3+119526x^2-42957x+2061.
$$
For lower bounds on $\gb$ one has to exhibit infinitely many
non-isomorphic curves over a fixed number field $K$, having  a given number
of  $K$-rational points. There are a few  different methods that have been
used for this purpose. One basic technique is to find a
surface (in
$\P ^3$ say) containing as many lines
as possible, and then construct a family of curves by taking plane sections
of that surface. The points of intersection with the lines are the seeked for rational points. 

This 
method has been
exploited in \cite{CHMjr} and by Noam Elkies, to give
 the (partly unpublished)  results listed below.

\

$N_2> 128$

$N_3> 100$

$N_4> 126$

$N_5> 132$

$N_6> 146$

$N_9> 180$

$N_{10}> 192$

$N_{45}> 781$

\
 
\ni
Different techniques yield results 
for all genera: the bound below was obtained independently by Brumer and
Mestre

$$\gb \ge 16(g+1)$$

\ni
A description of Mestre's method can be found in [CHMjr],
while that of Brumer is in [C1]. Although the methods of Brumer and Mestre 
look very
 different from each other,   Brumer  discovered a precise relation between them. 
Notice also that the above 
bound is not sharp, in fact Elkies got
$N_{45}\geq 781$ and, of course, $781 >16(45+1)$.

\

\subsection{A function field analog.} 

The uniformity questions 
that we so far discussed  only over number fields, remain interesting for curves over function fields.
The rest of the paper is devoted to this case, and we will from now on consider
curves defined over a complex function field $L=\C(B)$ where $B$ is a nonsingular projective variety.
The main conjecture to keep in mind is the following (discussed for number fields in the previous sections):

\begin{conj}
For any complex function field $L$
and any integer $g\geq 2$ there exists a number $B_g(L)$ such that for every $C\in \CL$ we have
$|C(L)|\leq B_g(L)$ 
\end{conj}

A revision of  Parshin's  Lemma~\ref{trick}, using the moduli theory of curves, yields that it actually holds
uniformly; in other words,  if a uniform  Shafarevich conjecture holds,
then a  uniform Mordell conjecture holds. 

Therefore a natural attempt is to focus on sets of curves over function fields,
looking for their uniformity properties.

The Shafarevich problem can be approached within the framework of moduli theory.
Given a curve $C$ of genus $g$ defined over $L$ (always considered modulo $L$-isomorphism),
pick a family $f:X\la B$ whose generic fiber is $L$-isomorphic to $C$.
We can associate to such a family its classifying map (see \ref{mod})
$$
\gamma _f : B \da M_g ;
$$
since  such a map is uniquely determined by the curve $C$
(if  $f':X'\la B'$ is another family having $C$ as generic fiber,
then $\gamma _{f'}$ coincides with $\gamma _f$ on some open subset of $B$), we shall also 
use the notation $\gamma
_C:=\gamma _f$. 

In this situation sets of maps are  easier to handle than sets of curves,
we shall therefore 
mimick Definition~\ref{curves} and set

\begin{defi}
\label{maps}  Let $L=\C (B)$ be a  complex function field.
\begin{enumerate}  \item Denote by $\ML$ the set  of rational maps
$\gamma :B\da M_g$ such that there exists a curve $C\in \CL$ having $\gamma$ as classifying map.

\item 
Let $S$ be a closed subset of $B$.
Denote by $\MLS$ the subset of $\ML$ parametrizing  classifying maps of curves in $\CL$  having good 
 reduction outside of $S$
\end{enumerate} 
\end{defi}

We view  $\ML$ as the set of $L$-rational points of the moduli functor
(or moduli stack) of smooth   curves of genus $g$.
 
There is a canonical surjection
$$
\begin{array}{lccr}
\mu : & \CL & \la & \; \ML\\
\; & C &\mapsto   &\gamma _C\\
\end{array}
$$
and $\MLS$ is the image via $\mu$ of the locus of curves in $\CL$ for which  there exists 
a model $f:X\la B$ having nonsingular  fiber over every point of $B$
which does not lie in $S$.

The crucial observation is that the map $\mu$  has finite fibers;
more precisely, 
there exists a universal constant $D(g)$ (independent of $L$)
such that the fibers of $\mu$ have cardinality at most $D(g)$.
This follows from  foundational moduli theory:
given $n$ different curves $C_1, \ldots, C_n $ all defined over $L$ and
having the same classifying map 
$$\gamma _{C_1} = \ldots =\gamma _{C_n} = \gamma :B \da M_g$$
there exists a finite 
extension $L'$ of $L$ over which the $n$ curves above become all isomorphic
(\cite{DM}). Since the curves  are distinct (by assumption) we find $n$ different $L'$ automorphisms
for a curve defined over $L'$.
If $n$ is too large, this gives a contradiction
with the fact that the number of automorphism of a curve of genus $g$  is bounded by a constant $D(g)$
that depends only on $g$ (and not on the field of definition).

The conclusion is that uniform boundedness  of sets of curves is equivalent
to uniform boundedness  of rational points of the moduli functor.
Combining this with lemma~\ref{trick} 
(which, as we said, holds uniformly)
we obtain that
uniform boundedness of rational points on curves would follow from uniform boundedness of
rational points of the moduli  functor.

The  general result that one obtains in this framework is the following (\cite{C3}):

\begin{thm}
\label{usm} Fix integers $g\geq 2$, $d$ and $s$.
There exist numbers $M_g(d,s)$ and $B_g(d,s)$ such that, for every projective variety $B\subset \P ^r$ of degree
at most $d$, and for every closed subset $S$ of $B$ having degree  at most $s$, we have
\begin{enumerate}
\item  $|{\mathcal M}_g(L, S)|\leq M_g(d,s)$ \  ($L=\C(B)$ as always).
\item 
 For every $C\in {\mathcal C}_g(L, S)$
we have $|C(L)|\leq B_g(d,s)$.
\end{enumerate}
\end{thm} 
The two bounds above do not depend on  the dimension of $B$ or of the ambient projective space.
\begin{proof}
We outline the argument.
The first step  is the case $\dim B=1$, which is the difficult part.
One proves that  there exists a bound $P_g(q,s)$ 
such that for every curve $B$ of genus at most $q$,
and for every subset $S$ of $B$ having at most $s$ points,
we have $|{\mathcal M}_g(\C(B), S)|\leq P_g(q,s)$.

This is  the crucial point, achieved by constructing a moduli space 
over $M_g$ for the maps parametrized by all
sets $\MLS$ (see \cite{C2}).

As we explained before, this gives (is equivalent to) an identical
statement for  ${\mathcal C}_g(L, S)$.
Then one applies the uniform version of the Kodaira-Parshin trick to obtain
the proof of part 2 of the theorem for $\dim B =1$.

The higher dimensional statements are both obtained by ``slicing" $B$ by hyperplanes to obtain one-dimensional
sections. To such curves (whose genus is obviously bounded above by a function of $d$ only)
one applies the one-dimensional result to conclude the proof.
\end{proof}

\begin{remark} 
In a very different framework, 
Miyaoka  (in \cite{Mi}) 
used the theory of vector bundles on surfaces
to obtain effective results of the same type as part 2,
for rational points of curves defined over one-dimensional function
fields and satisfying some other geometric assumptions.

An approach to the Shafarevich conjecture, using Chow varieties rather than
moduli spaces, has been recently applied by
Heier (\cite{He}) in his Ph.D. thesis, to produce an effective expression
for the bound $P_g(q,s)$ above (on the cardinality
of ${\mathcal C}_g(\C(B), S)$ when $B$ has dimension $1$ and genus $q$).
\end{remark}

\

It is at this point natural (at least to my opinion) to expect analogous 
boundedness phenomena to
hold in the arithmetic case; for instance, a statement like the following is not known:
{\it there exists a number $M_g(d,s)$ such that for every number field $K$ with $[K:\Q]\leq d$,
for every set $S$ of places  of $K$, with $|S|\leq s$, for every curve 
$C\in {\mathcal C}_g(K,S)$,
we have $|C(K)|\leq M_g(d,s)$. }

\

We conclude with  a few  very special cases in which the Conjecture holds. 
Let $L=\C (B)$ be our  complex function field.
\begin{enumerate}  \item 
Denote by ${\mathcal C}^2_g(L)$ the set of curves in 
$\CL$ 
having good reduction in codimension $1$,
that is 
$$
{\mathcal C}^2_g(L) := \bigcup _{codim _B S \geq 2} \CLS
$$
\item 
Denote by ${\mathcal C}^{max}_g(L)$ the set of curves $C$ in 
$\CL$ 
having maximal variation of moduli,
that is, the associated moduli map $\gamma _C$ is generically finite. 
\end{enumerate}

Then we have
\begin{prop}
\begin{enumerate}
\item There exists a number $B_g^2(L)$ such that for every curve $C \in {\mathcal C}^2_g(L)$
we have $|C(L)|\leq B_g^2(L)$.
\item
Fix $g\geq 24$ and let $trdeg _{\C}L =3g-3$. There exists a number $B_g^{max}(L)$ such that for every curve 
$C \in
{\mathcal C}^{max}_g(L)$ we have $|C(L)|\leq B_g^{max}(L)$.
\end{enumerate} 
\end{prop}
Both statements follow trivially from the fact that the sets of curves in question are finite.
The finiteness of ${\mathcal C}^2_g(L)$ can be proved in the same spirit as Theorem~\ref{usm}
(see \cite{C3}). 

The finiteness of ${\mathcal C}^{max}_g(L)$ is obtained using the well known theorem
of Harris-Mumford and Eisenbud-Harris
that $\Mg$ 
is of general type for $g\geq 23$ (see \cite{HM} for proof and references).
Combining such a fact with theorem \ref{KO} we obtain that the set of rational dominant maps
from $B$ to $M_g$ is finite. Now,  ${\mathcal C}^{max}_g(L)$ has a natural surjection
(the restriction of $\mu$) onto such a set, hence it is finite.

\

\noindent
{\bf{Aknowledgements.}}
It is a pleasure to thank Fr\'ed\'eric Campana and Burt Totaro for their kind comments and
suggestions.

\bigskip

\noindent Lucia Caporaso \\
Dipartimento di Matematica, Universit\`a  Roma Tre \\
Largo S.\ L.\ Murialdo, 1 \\
00146 Roma - ITALY\\
caporaso@mat.uniroma3.it

\
\end{document}